\documentclass[11pt, leqno]{article}
\usepackage{euscript,amstext,amsmath,amsthm,a4,amssymb}
\numberwithin{equation}{section}

\def\bC{{\rm \bf C}}

\def\bZ{{\rm \bf Z}}

\newcommand{\tr}{{\rm Tr}}

\newcommand{\comment}[1]{}

\newtheorem{thm}{Theorem}[section]
\newtheorem{prop}[thm]{Proposition}
 \newtheorem{cor}[thm]{Corollary}

\theoremstyle{remark}
\newtheorem{Rem}[thm]{Remark}
\theoremstyle{definition}

\newtheorem{ex}[thm]{Example}

\title{$\eta$-invariant  and flat vector bundles II}
\author{Xiaonan Ma\footnote{Centre de Math\'ematiques Laurent Schwartz,
UMR 7640 du CNRS,
\'Ecole Polytechnique, 91128 Palaiseau Cedex, France.
(ma@math.polytechnique.fr)}\ \ and\ \ Weiping Zhang\footnote{Chern
Institute of Mathematics \& LPMC, Nankai University, Tianjin
300071, P.R. China. (weiping@nankai.edu.cn)}}
\date{
}

\begin{document}

\maketitle
\begin{abstract} We first apply the method and results in the
previous paper to give a new proof of a result (hold in $ {\bf
C}/{\bf Z}$) of Gilkey on the
 variation of $\eta$-invariants
associated to non self-adjoint Dirac type operators. We then give
an explicit local expression of certain $\eta$-invariant appearing
in recent papers of Braverman-Kappeler on what they call refined
analytic torsion, and propose an alternate formulation of their
definition of the refined analytic torsion. A refinement in ${\bf
C}$ of the above variation formula is also proposed.

$\ $

{\bf Keywords} flat vector bundle, $\eta$-invariant, refined
analytic torsion

{\bf 2000 MR Subject Classification} 58J
\end{abstract}

\renewcommand{\theequation}{\thesection.\arabic{equation}}
\setcounter{equation}{0}

\section{Introduction} \label{s1}

In a previous paper \cite{MZ}, we have given an alternate
formulation of (the mod {\bf Z} part of) the $\eta$-invariant of
Atiyah-Patodi-Singer \cite{APS1, APS2, APS3} associated to
non-unitary flat vector bundles by identifying explicitly its real
and imaginary parts.

On the other hand,  Gilkey
  has studied this kind of $\eta$-invariants systematically in
  \cite{G}, and in particular proved a general variation formula
  for them. However, it lacks in \cite{G} the identification of the real and imaginary parts
  of the $\eta$-invariants as we did in \cite{MZ}.

In this article, we   first  show that our results in \cite{MZ}
lead to a direct derivation of Gilkey's variation formula
\cite[Theorem 3.7]{G}.

 The second purpose of this paper is to apply the results
in \cite{MZ} to examine the $\eta$-invariants appearing in the
recent papers of Braverman-Kappeler \cite{BK1, BK2, BK3} on
refined analytic torsions. We show that the imaginary part of the
$\eta$-invariant appeared in these articles admits an explicit
local expression which suggests an alternate formulation of the
definition of the refined analytic torsion there. This
reformulation provides an analytic resolution of a problem due to
Burghelea (\cite{BuH0, BuH}) on the existence of a univalent
holomorphic function on the representation space having the
Ray-Singer analytic torsion as its absolute value.

Finally, using the extension  (to the case of non-self-adjoint
operators) given in \cite{ZL} of the concept of spectral flow
\cite{APS3}, we propose a refinement in ${\bf C}$ of the above
variation formula for $\eta$-invariants.

$\  $

\noindent {\bf Acknowledgements}. We would like to thank Maxim
Braverman for bringing \cite{G} to our attention. The work of the
second author was partially supported by the National Natural
Science Foundation of China.

\section{$\eta$-invariant and the variation  formula}\label{s2}

Let $M$ be an odd dimensional oriented closed  spin manifold
carrying a Riemannian metric $g^{TM}$. Let $S(TM)$ be the
associated Hermitian vector bundle of spinors. Let $(E, g^E)$ be a
Hermitian vector bundle over $M$ carrying a unitary connection
$\nabla^E$. Moreover, let $(F,g^F)$ be a Hermitian vector bundle
over $M$ carrying a
  flat connection $\nabla^F$. We do not assume that $\nabla^F$
preserves the Hermitian metric $g^F$ on $F$.

   Let
$ D^{E\otimes F}:\Gamma(S(TM)\otimes E\otimes F)\longrightarrow
\Gamma(S(TM)\otimes E\otimes F)$
denote the corresponding (twisted) Dirac operator.

It is pointed  out in \cite[Page 93]{APS3} that one can define the
reduced $\eta$-invariant of  $D^{E\otimes F}$, denoted by
$\overline{\eta}( D^{E\otimes F})$,  by working on (possibly)
non-self-adjoint elliptic operators.

In this section, we will first recall the main result in \cite{MZ}
on $\overline{\eta}( D^{E\otimes F})$ and then show how it leads
directly to a proof of the variation formula of Gilkey \cite[Theorem 3.7]{G}.

\subsection{Chern-Simons classes and flat vector
bundles}\label{2a} \setcounter{equation}{0} We fix a square root
of $\sqrt{-1}$ and let
$\varphi:\Lambda(T^*M)\rightarrow\Lambda(T^*M)$ be the
homomorphism defined by
$\varphi:\omega\in\Lambda^i(T^*M)\rightarrow (2\pi
\sqrt{-1})^{-i/2}\omega.$ The formulas in what follows will not
depend on the choice of the square root of $\sqrt{-1}$.

If  $W$ is a  complex vector bundle over $M$ and  $\nabla^W_0$,
$\nabla^W_1$ are two connections on $W$. Let $W_t$, $0\leq t\leq
1$, be a smooth path of connections on $W$ connecting $\nabla^W_0$
and $\nabla^W_1$. We define the Chern-Simons form
$CS(\nabla^W_0,\nabla^W_1)$ to be the differential form given by
\begin{align}\label{2.5}
CS\left(\nabla^W_0,\nabla^W_1\right)=-\left(1\over
2\pi\sqrt{-1}\right)^{1\over 2}
\varphi\int_0^1\tr\left[{\partial\nabla^W_t\over\partial
t}\exp\left(-\left(\nabla^W_t\right)^2\right) \right]dt.
\end{align}
Then (cf. \cite[Chapter 1]{Z1})
\begin{align}\label{2.6}
dCS\left(\nabla^W_0,\nabla^W_1\right)={\rm
ch}\left(W,\nabla^W_1\right)-{\rm ch}\left(W,\nabla^W_0\right).
\end{align}
Moreover, it is well-known that up to exact forms,
$CS(\nabla^W_0,\nabla^W_1)$ does not depend on the path of
connections on $W$ connecting $\nabla^W_0$ and $\nabla^W_1$.

 Let $(F,\nabla^F)$ be
a flat vector bundle carrying the flat connection $\nabla^F$. Let
$g^F$ be a Hermitian metric on $F$. We do not assume that
$\nabla^F$ preserves $g^F$. Let $(\nabla^F)^*$ be the adjoint
connection of $\nabla^F$ with respect to $g^F$.

From \cite[(4.1), (4.2)]{BZ} and \cite[\S 1(g)]{BL}, one has
\begin{align}\label{2.1}
\left(\nabla^F\right)^*=\nabla^F+\omega\left(F,g^F\right)
\end{align}
with
\begin{align}\label{2.2}
 \omega\left(F,g^F\right)=\left(g^F\right)^{-1}\left(\nabla^Fg^F\right).
\end{align} Then
\begin{align}\label{2.3}
  \nabla^{F,e}=\nabla^F+{1\over 2}\omega\left(F,g^F\right)
\end{align}
is a Hermitian connection on $(F,g^F)$ (cf. \cite[(4.3)]{BZ}).

Following \cite[(2.6)]{MZ} and \cite[(2.47)]{MZ1}, for any
$r\in{\bf C}$, set
\begin{align}\label{2.4}
  \nabla^{F,e,(r)}=\nabla^{F,e}+{\sqrt{-1}r\over
  2}\omega\left(F,g^F\right).
\end{align}
Then for any $r\in{\bf R}$, $\nabla^{F,e,(r)}$ is a Hermitian
connection on $(F,g^F)$.

On the other hand, following \cite[(0.2)]{BL},  for any integer
$j\geq 0$, let $c_{2j+1}(F,g^F)$ be the Chern form defined by
\begin{align}\label{2.7}
c_{2j+1}\left(F,g^F\right)=\left(2\pi\sqrt{-1}\right)^{-j}
2^{-(2j+1)}\tr\left[\omega^{2j+1}\left(F,g^F\right)\right].
\end{align} Then $c_{2j+1}(F,g^F)$ is a closed form on $M$.
Let $c_{2j+1}(F)$ be the associated cohomology class in
$H^{2j+1}(M,{\bf R})$, which does not depend on the choice of
$g^F$.

For any $j\geq 0$ and $r\in {\bf R}$, let $a_j(r)\in {\bf R}$ be
defined as
\begin{align}\label{2.8}
a_j(r)=\int_0^1\left(1+u^2r^2\right)^jdu.
\end{align}

With these notation we can now state the following result first
proved in \cite[Lemma 2.12]{MZ1}.

\begin{prop}\label{t2.1} The following identity in
$H^{\rm odd}(M,{\bf R})$ holds for any  $r\in {\bf R}$,
\begin{align}\label{2.9}
CS\left(\nabla^{F,e},\nabla^{F,e,(r)}\right)=-{r\over
2\pi}\sum_{j=0}^{+\infty} {a_j(r)\over j!}c_{2j+1}(F).
\end{align}
\end{prop}

\subsection{$\eta$   invariant associated to flat vector
bundles}\label{2b}

Let
\begin{align}\label{2.10}D^{E\otimes F,e} :\Gamma(S(TM)\otimes E\otimes F)\longrightarrow
\Gamma(S(TM)\otimes E\otimes F)
\end{align}
 denote the Dirac operator associated to the  connection
 $\nabla^{F,e }$ on $F$ and $\nabla^E$ on $E$. Then $D^{E\otimes
 F,e} $ is formally self-adjoint and one can define the associated
 reduced $\eta$-invariant as in \cite{APS1}.

 In view of Proposition \ref{t2.1}, one can restate the main
 result of \cite{MZ}, which is \cite[Theorem 2.2]{MZ}, as follows,
 \begin{align}\label{2.14}
\overline{\eta}\left(D^{E\otimes F}\right)\equiv
\overline{\eta}\left(D^{E\otimes F,e}\right)+
\int_M\widehat{A}(TM){\rm ch}(E)CS\left(\nabla^{F,e},\nabla^{F
}\right)\ \ {\rm mod}\ {\bf Z},
\end{align}
where $\widehat{A}(TM)$ and ${\rm ch}(E)$ are the $\widehat{A}$
class of $TM$ and the Chern character of $E$ respectively (cf.
\cite{Z1}).

Now let $\widetilde{\nabla}^F$ be another flat connection on $F$.
We use the notation with $\widetilde{\ }$ to denote the objects
associated with this flat connection.

Then one has
\begin{align}\label{2.15}
\overline{\eta}\left(\widetilde{D}^{E\otimes F}\right)\equiv
\overline{\eta}\left(\widetilde{D}^{E\otimes F,e}\right)+
\int_M\widehat{A}(TM){\rm
ch}(E)CS\left(\widetilde{\nabla}^{F,e},\widetilde{\nabla}^{F
}\right)\ \ {\rm mod}\ {\bf Z}.
\end{align}

By the variation formula for $\eta$-invariants associated to
self-adjoint Dirac operators (\cite{APS1}, \cite{BF}), one knows
that
\begin{align}\label{2.16}
\overline{\eta}\left(\widetilde{D}^{E\otimes
F,e}\right)-\overline{\eta}\left(D^{E\otimes
F,e}\right)\equiv\int_M\widehat{A}(TM){\rm
ch}(E)CS\left(\nabla^{F,e},\widetilde{\nabla}^{F,e}\right)\ \ {\rm
mod}\ {\bf Z}.
\end{align}

From (\ref{2.14})-(\ref{2.16}), one deduces that
\begin{align}\label{2.17}
\overline{\eta}\left(\widetilde{D}^{E\otimes
F}\right)-\overline{\eta}\left(D^{E\otimes F}\right)
\equiv\int_M\widehat{A}(TM){\rm
ch}(E)CS\left(\nabla^{F,e},\widetilde{\nabla}^{F,e}\right)\end{align}
$$-\int_M\widehat{A}(TM){\rm
ch}(E)CS\left(\nabla^{F,e},\nabla^{F }\right)+
\int_M\widehat{A}(TM){\rm
ch}(E)CS\left(\widetilde{\nabla}^{F,e},\widetilde{\nabla}^{F
}\right)$$
 $$=\int_M\widehat{A}(TM){\rm ch}(E)CS\left(\nabla^{F
},\widetilde{\nabla}^{F } \right)\ \ {\rm mod}\ {\bf Z},$$ which
is exactly the Gilkey formula \cite[Theorem 1.6]{G}
for the operator $P=D^E$ therein.

\begin{Rem}\label{100} As was indicated in \cite[Remark 2.4]{MZ}, the main
result in \cite{MZ} holds also for general Hermitian vector
bundles equipped with a (possibly) non-Hermitian connection.
Indeed, if we do not assume that $\nabla^F$ is flat, then at
least (\ref{2.1})-(\ref{2.4}) still holds. Thus for any $r\in{\bf
R}$, we have well-defined (formally self-adjoint) operator
$D^{E\otimes F}(r)$ which is associated to the Hermitian
connection $\nabla^{F,e,(r)}$ on $F$. For any $r\in{\bf R}$, one
then has the variation formula (cf. \cite{APS1} and \cite{BF})
\begin{align}\label{2.18}
\overline{\eta}\left( {D}^{E\otimes
F}(r)\right)-\overline{\eta}\left(D^{E\otimes F,e}\right)
\equiv\int_M\widehat{A}(TM){\rm ch}(E)CS\left(\nabla^{F,e},
{\nabla}^{F,e,(r)}\right)\ \ {\rm mod}\ {\bf Z}.\end{align} By
(\ref{2.5}), one sees easily that the right hand side of
(\ref{2.18}) is a holomorphic function (indeed a polynomial) of
$r$. Thus, by analytic continuity, as in \cite{MZ}, one gets that
for any $r\in {\bf C}$, (\ref{2.18}) still holds. In particular,
if we set $r=\sqrt{-1}$, we get
\begin{align}\label{2.28}
\overline{\eta}\left( {D}^{E\otimes F}
\right)\equiv\overline{\eta}\left(D^{E\otimes F,e}\right)
+\int_M\widehat{A}(TM){\rm ch}(E)CS\left(\nabla^{F,e}, {\nabla}^{F
}\right)\ \ {\rm mod}\ {\bf Z},\end{align} which generalizes
(\ref{2.14}). Then by proceeding as above, we see that
(\ref{2.17}) holds without the assumption of the flatness of
connections $\nabla^F$ and $\widetilde{\nabla}^F$.

By \eqref{2.5} and  \eqref{2.4},
\begin{multline}\label{a2.20}
CS\left(\nabla^{F,e}, {\nabla}^{F,e,(r)}\right)=
\frac{-1}{2\pi}
\int_0^1 \tr \left[\frac{r}{2}\omega\left(F,g^F\right)
\exp\left(\frac{-1}{2\pi\sqrt{-1}}\left(\nabla^{F,e,(tr)}\right)^2\right)\right]
dt\\
=\sum_{i=0}^{\dim M} a_i\left(\nabla^F, g^F\right)\, r^i.
 \end{multline}

By (\ref{2.4}), one has
\begin{align}\label{2.19}
\left(\nabla^{F,e,(r)}\right)^2=\left(\nabla^{F,e }\right)^2 +
\frac{\sqrt{-1}r}{2}\left(\nabla^{F,e
}\omega\left(F,g^F\right)\right)-\frac{r^2}{4}\left(\omega\left(F,g^F\right)\right)^2
 .\end{align}
Note that
\begin{align}\label{a2.19}
\nabla^{F,e}\omega\left(F,g^F\right)
= [\nabla^{F,e}, \omega\left(F,g^F\right)]=0, \quad {\rm if} \,\,
\nabla^F \,\, \,\, \mbox{is flat}.
\end{align}

 By taking adjoint of \eqref{2.19},
we see that when $r\in{\bf C}$ is purely imaginary, one has
\begin{multline}\label{2.20}
 \left({1\over 2\pi\sqrt{-1}}
\left(\nabla^{F,e,(r)}\right)^2\right)^*= \frac{1}{2\pi\sqrt{-1}}\left(\left(\nabla^{F,e }\right)^2 -
\frac{r^2}{4}\left(\omega\left(F,g^F\right)\right)^2\right)\\
-\left(\frac{1}{2\pi\sqrt{-1}}\right) \frac{\sqrt{-1}r}{2}\left(\nabla^{F,e}
\omega\left(F,g^F\right)\right).
  \end{multline}

From (\ref{2.5}), (\ref{a2.20})
and (\ref{2.20}), one sees that when $r\in{\bf C}$ is purely
imaginary, then
\begin{align}\label{2.21}
\begin{split}
& {\rm
Re}\left(CS\left(\nabla^{F,e},{\nabla}^{F,e,(r)}\right)\right)
=\sum_{i\ {\rm even}} a_i\left(\nabla^F, g^F\right)\, r^i,\\
&{\rm
Im}\left(CS\left(\nabla^{F,e},{\nabla}^{F,e,(r)}\right)\right)
=\frac{1}{\sqrt{-1}}\sum_{i\ {\rm odd}} a_i\left(\nabla^F,
g^F\right)\, r^i.
 \end{split} \end{align}
Thus  when $r\in{\bf C}$ is purely imaginary,
from \eqref{2.28} and \eqref{2.21}, we have
\begin{align}\label{2.22}
\begin{split}
& {\rm Re}\left(\overline{\eta}\left({D}^{E\otimes
F}(r)\right)\right)\equiv\overline{\eta}\left(D^{E\otimes
F,e}\right)+\sum_{i\ {\rm even}} r^i
\int_M\widehat{A}(TM){\rm ch}(E)a_i\left(\nabla^F, g^F\right)\ \ {\rm mod}\ {\bf Z},\\
& {\rm Im}\left(\overline{\eta}\left({D}^{E\otimes
F}(r)\right)\right)= \frac{1}{\sqrt{-1}}\sum_{i\ {\rm
odd}}r^i\int_M\widehat{A}(TM){\rm ch}(E)a_i\left(\nabla^F,
g^F\right).
 \end{split} \end{align}

In particular, by setting $r=\sqrt{-1}$, we get
\begin{align}\label{2.23}
\begin{split}
& {\rm Re}\left(\overline{\eta}\left({D}^{E\otimes F} \right)\right)
\equiv\overline{\eta}\left(D^{E\otimes F,e}\right)
+\sum_{i\ {\rm even}}(-1)^{i\over 2}
\int_M\widehat{A}(TM){\rm ch}(E)a_i\left(\nabla^F, g^F\right) \ \ {\rm mod}\ {\bf Z},\\
&{\rm Im}\left(\overline{\eta}\left({D}^{E\otimes F}\right)\right)
= \sum_{i\ {\rm odd}}\left( {-1}\right)^{i-1\over
2}\int_M\widehat{A}(TM){\rm ch}(E)a_i\left(\nabla^F, g^F\right).
 \end{split} \end{align}

 This generalizes the main result
in \cite{MZ}.
\end{Rem}

\section{$\eta$-invariant and the refined analytic torsion of\\
Braverman-Kappeler}\label{s3}

Recently, in a series of preprints \cite{BK1, BK2, BK3}, Braverman
and Kappeler introduce what they call refined analytic torsion.
The $\eta$-invariant associated with flat vector bundles plays a
role in their definition. In this section, we first examine the
imaginary part of the $\eta$-invariant appearing in \cite{BK1,
BK2, BK3} from the point of view of the previous sections and
propose an alternate definition of the refined analytic torsion.
We then combine this refined analytic torsion with the
$\eta$-invariant to construct analytically a univalent holomorphic
function on the space of   representations of $\pi_1(M)$ having
the absolute value equals to the Ray-Singer torsion, thus
resolving a problem posed by Burghelea (cf. \cite{BuH}).

\subsection{$\eta$-invariant and the refined analytic torsion of
Braverman-Kappeler}\label{s3.1}

Since there needs no spin condition in \cite{BK1, BK2, BK3}, here
we   start with a closed  oriented smooth  odd dimensional
manifold $M$ with $\dim M=2n+1$. Let $g^{TM}$ be a Riemannian
metric on $TM$. For any $X\in TM$, let $X^*\in T^*M$ denote its
metric dual and $c(X)=X^*-i_X$ denote the associated Clifford
action acting on $\Lambda^*(T^*M)$, where $X^*$ and $i_X$ are the
notation for the exterior and interior multiplications  of $X$
respectively.

Let $e_1$, $\dots$, $e_{2n+1}$ be an oriented orthonormal basis of
$TM$.  Set
\begin{align}\label{3.1}
\Gamma=\left(\sqrt{-1}\right)^{n+1 }c(e_1)\cdots c(e_{2n+1}).
\end{align}
Then $\Gamma^2={\rm Id}$ on $\Lambda^*(T^*M)$.

Let $(F,g^F)$ be a Hermitian vector bundle over $M$ equipped with a flat
connection $\nabla^F$ which need not preserve the Hermitian metric
$g^F$ on $F$. Then the exterior differential $d$ on
$\Omega^*(M)=\Gamma(\Lambda^*(T^*M))$ extends naturally to the
twisted exterior differential $d^F$ acting on
$\Omega^*(M,F)=\Gamma(\Lambda^*(T^*M)\otimes F)$.

We define the twisted signature operator $D_{\rm Sig}^F$ to be
\begin{align}\label{3.2}
D_{\rm Sig}^F={1\over 2}\left(\Gamma
d^F+d^F\Gamma\right):\Omega^{\rm even}(M,F)\rightarrow \Omega^{\rm
even}(M,F).
\end{align}
It coincides with the odd signature operator ${1\over 2}{\cal
B}_{\rm even} $ in \cite{BK1, BK2, BK3}.

Let $\nabla^{\Lambda^{\rm even}(T^*M)\otimes F}$
(resp. $\nabla^{\Lambda^{\rm even}(T^*M)\otimes F,e}$) be the tensor
product connections on $\Lambda^{\rm even}(T^*M)\otimes F$ obtained
from $\nabla^F$ (resp. $\nabla^{F,e}$) and the canonical connection on
$\Lambda^{\rm even}(T^*M)$ induced by the Levi-Civita connection
$\nabla^{TM}$ of $g^{TM}$.

From (\ref{3.2}), it is easy to verify that
\begin{align}\label{3.3}
D_{\rm Sig}^F
=\Gamma\left(\sum_{i=1}^{2n+1}c(e_i)\nabla^{\Lambda^{\rm
even}(T^*M)\otimes F}_{e_i}\right).
\end{align}

Set
\begin{align}\label{3.3a}
D_{\rm Sig}^{F,e}
=\Gamma\left(\sum_{i=1}^{2n+1}c(e_i)\nabla^{\Lambda^{\rm
even}(T^*M)\otimes F,e}_{e_i}\right).
\end{align}
Then $D_{\rm Sig}^{F,e}$ is formally self-adjoint.

Since locally one has identification $S(TM)\otimes
S(TM)=\Lambda^{\rm even}(T^*M)$, one sees that one can apply the
results in the previous section to the case $E=S(TM)$ to the
current situation.

In particular, we get
\begin{align}\label{3.4}
\begin{split}
{\rm Re}\left(\overline{\eta}\left(D_{\rm Sig}^F\right)\right)
&\equiv \overline{\eta}\left(D_{\rm Sig}^{F,e}\right)
 \quad\ \ \ \ \ \ \ \ \   \mbox{mod}\ \bZ,\\
{\rm Im}\left(\overline{\eta}\left(D_{\rm Sig}^F\right)\right)
&={1\over \sqrt{-1}}\int_M {\rm
L}\left(TM,\nabla^{TM}\right)CS\left(\nabla^{F,e},\nabla^{F }\right)\\
&= -{1\over 2\pi}\int_M{\rm
L}(TM )\sum_{j=0}^{+\infty}{2^{2j}j!\over (2j+1)!}c_{2j+1}(F),
\end{split}
\end{align}
where ${\rm L}(TM,\nabla^TM)$ is the Hirzebruch ${\rm L}$-form
defined by
\begin{align}\label{3.5}
{\rm L}\left(TM,\nabla^TM\right)=\varphi\, {\rm
det}^{1/2}\left({R^{TM}\over \tanh \left(R^{TM}/2\right)}\right),
\end{align}
with $R^{TM}=(\nabla^{TM})^2$ the curvature of $\nabla^{TM}$, and
${\rm L}(TM)$ is the associated class.

$\ $

\begin{Rem}
By proceeding as in Section 2, we can get \cite[Theorem 3.7]{G}
easily by using the results in Remark \ref{100}.
\end{Rem}



\begin{prop}\label{t3.1} The function
\begin{align}\label{3.6}
\Psi\left(F,\nabla^F\right)={\rm
Im}\left(\overline{\eta}\left(D_{\rm Sig}^F\right)\right) +{1\over
2\pi}\int_M{\rm L}(TM ) c_{ 1}(F)
\end{align}
is a locally constant function on the set of flat connections on
$F$. In particular, $\Psi (F,\nabla^F )=0$ if $\nabla^F$ can be
connected to a unitary flat connection through a path of flat
connections.
\end{prop}

{\it Proof}. Let $\nabla_t^F$, $0\leq t\leq 1$, be a smooth pass
of flat connections on $F$.

From (\ref{3.4}), we get
\begin{align}\label{3.7}
 { \sqrt{-1}}{\rm Im}\left(\overline{\eta}\left(D_{\rm Sig,1}^F\right)\right)
 -{ \sqrt{-1}}{\rm Im}\left(\overline{\eta}\left(D_{\rm Sig,0}^F\right)\right)
\end{align}
$$=\int_M {\rm L}\left(TM,\nabla^{TM}\right)CS\left(\nabla^{F,e}_1,\nabla^{F
}_1\right)-\int_M {\rm
L}\left(TM,\nabla^{TM}\right)CS\left(\nabla^{F,e}_0,\nabla^{F
}_0\right)$$
$$=\sqrt{-1}\int_M {\rm L}\left(TM,\nabla^{TM}\right){\rm Im}\left(CS\left(\nabla^{F,e}_1,\nabla^{F,e}_0
\right)- CS\left(\nabla^{F }_1,\nabla^{F }_0\right)\right)$$
$$=\sqrt{-1}\int_M {\rm L}\left(TM,\nabla^{TM}\right){\rm Im}\left(
 CS\left(\nabla^{F }_0,\nabla^{F }_1\right)\right).$$

 Now consider the path of flat connections $\nabla_t^F$, $0\leq t\leq 1$.
 Since for any $t\in [0,1]$, $(\nabla_t^F)^2=0$, from (\ref{2.5}),
 (\ref{2.3}),
 one gets
\begin{align}\label{3.8}
 CS\left(\nabla^{F }_0,\nabla^{F }_1\right)=\left({1\over
 2\pi\sqrt{-1}}\right)\left(\nabla^F_0-\nabla^F_1\right)
\end{align}
$$=\left({1\over
 2\pi\sqrt{-1}}\right)\left(\nabla^{F,e}_0-\nabla^{F,e}_1\right)
 -\left({1\over
 2\pi\sqrt{-1}}\right)\left({1\over 2}\omega_0(F,g^F)-{1\over
 2}\omega_1(F,g^F)\right).$$

Thus, one has
\begin{align}\label{3.9}
\sqrt{-1}{\rm Im}\left(
 CS\left(\nabla^{F }_0,\nabla^{F }_1\right)\right)=-{1\over 2\pi\sqrt{-1}}\left({1\over 2}\omega_0(F,g^F)-{1\over
 2}\omega_1(F,g^F)\right)
\end{align}
$$=-{1\over 2\pi\sqrt{-1}}\left(c_1\left(F,\nabla^F_0\right)-c_1\left(F,\nabla^F_1\right)\right).$$

From (\ref{3.7}) and (\ref{3.9}), we get
\begin{align}\label{3.10}
{\rm Im}\left(\overline{\eta}\left(D_{\rm
Sig,1}^F\right)\right)+{1\over 2\pi }\int_M {\rm
L}\left(TM,\nabla^{TM}\right) c_1\left(F,\nabla^F_1\right)
\end{align}
$$={\rm
Im}\left(\overline{\eta}\left(D_{\rm
Sig,0}^F\right)\right)+{1\over 2\pi }\int_M {\rm
L}\left(TM,\nabla^{TM}\right) c_1\left(F,\nabla^F_0\right),$$ from
which Proposition \ref{t3.1} follows. \ \ Q.E.D.

$\ $

\begin{Rem}\label{r3.1} Formula (\ref{3.10}) is closely related to
\cite[Theorem 12.3]{BK1}. Moreover, for any representation
$\alpha$ of the fundamental group $\pi_1(M)$, let
$(F_\alpha,\nabla^{F_\alpha})$ be the associated flat vector
bundle. One has
\begin{align}\label{3.11}
 \exp\left(\pi\Psi\left(F_\alpha,\nabla^{F_\alpha}\right)\right)=r(\alpha),
\end{align}
where $r(\alpha)$ is the function appearing in \cite[Lemma
5.5]{BK3}. While from (\ref{3.4}) and (\ref{3.6}), one has
\begin{align}\label{3.12}
\Psi\left(F,\nabla^F\right)= -{1\over 2\pi}\int_M{\rm L}(TM
)\sum_{j=1}^{+\infty}{2^{2j}j!\over (2j+1)!}c_{2j+1}(F).
\end{align}
Combining with (\ref{3.11}), this gives an explicit local
expression of $r(\alpha)$ as  well as the locally constant function
$r_{\cal C}$ defined in \cite[Definition 5.6]{BK3}.
\end{Rem}

$\ $

\begin{Rem}\label{r3.2} To conclude this subsection, we
propose to modify the definition of the refined analytic torsion
of Braverman-Kappeler as follows: for any Hermitian vector bundle
equipped with a flat connection $\nabla^F$ over an  oriented
closed smooth odd dimensional manifold $M$ equipped with a Riemannian metric
$g^{TM}$, let $\rho(\nabla^F,g^{TM})$ be the element defined in
\cite[(2.13)]{BK3}. Then we propose the definition of the refined
analytic torsion as
\begin{align}\label{3.13}
 \rho_{\rm
 an}'\left(\nabla^F,g^{TM}\right)=
 \rho\left(\nabla^F,g^{TM}\right)e^{
 \pi\sqrt{-1}{\rm rk}(F)\overline{\eta}(D_{\rm sig})},
 \end{align} where $\overline{\eta}(D_{\rm sig})$ is the reduced
 $\eta$ invariant  in the sense of Atiyah-Patodi-Singer \cite{APS1}
 of the signature operator coupled with the trivial complex line
 bundle over $M$ (i.e. $D_{\rm sig}:=D_{\rm sig}^\bC$).
 \comment{There are two advantages of this reformulation. First, by
 multiplying the local factor
 $e^{-\pi\Psi\left(F,\nabla^F\right)}$ makes the comparison formula
\cite[(5.8)]{BK3} of
 the refined analytic torsion has closer resemblance  in comparing with
 the formulas of Cheeger-M\"uller and Bismut-Zhang (cf.
 \cite{BZ}).}
The advantage of this reformulation is that since
$\overline{\eta}(D_{\rm sig})$ various smoothly with respect to
the metric $g^{TM}$ (as the dimension of ${\rm ker}( D_{\rm sig})$
does not depend on the metric $g^{TM}$), the ambiguity of the
power of $\sqrt{-1}$ disappears if one uses $e^{\pi\sqrt{-1}{\rm
rk}(F)\overline{\eta}(D_{\rm sig})}$ to replace the factor
$e^{{\pi\sqrt{-1} {\rm rk}
 (F)\over 2}\int_N L(p,g^M)}$ in \cite[(2.14)]{BK3}.
\end{Rem}

\subsection{Ray-Singer analytic torsion and  univalent holomorphic functions on the
representation space}\label{s.3.2}

Let $(F,\nabla^F)$ be a complex flat vector bundle. Let $g^F$ be
an Hermitian metric on $F$. We fix a flat connection
$\widetilde{\nabla}^F$ on $F$ (Note here that we do not assume
that $\nabla^F$ and $\widetilde{\nabla}^F$ can be connected by a
smooth path of flat connections).

Let $g^{TM}$ be a Riemannian metric on $TM$ and $\nabla^{TM}$ be
the associated Levi-Civita connection.

Let $\widetilde{\eta}(\nabla^F,\widetilde{\nabla}^F)\in{\bf C}$ be
defined by
\begin{align}\label{3.16}
\widetilde{\eta}\left(\nabla^F,\widetilde{\nabla}^F\right)=\int_M{\rm
L}\left(TM,\nabla^{TM}\right)CS\left(\widetilde{\nabla}^{F,e},\nabla^F\right).
 \end{align}
 One verifies easily that $\widetilde{\eta}(\nabla^F,\widetilde{\nabla}^F)\in{\bf
 C}$ does not depend on $g^{TM}$, and is a holomorphic function of
 $\nabla^F$.
Moreover, by (\ref{3.4}) one has
\begin{align}\label{3.161}
{\rm
Im}\left(\widetilde{\eta}\left(\nabla^F,\widetilde{\nabla}^F\right)\right)=
{\rm Im}\left(\overline{\eta}\left(D_{\rm Sig}^F\right)\right).
 \end{align}

 Recall that we have modified the refined analytic torsion of
 \cite{BK1, BK2, BK3} in (\ref{3.13}).

 Set
\begin{align}\label{3.17}
 {\cal T}_{\rm
 an}\left(\nabla^F,g^{TM}\right)=
 \rho_{\rm an}'\left(\nabla^F,g^{TM}\right)
 \exp\left(\sqrt{-1}\pi\widetilde{\eta}\left(\nabla^F,\widetilde{\nabla}^F\right)\right).
 \end{align}
 Then ${\cal T}_{\rm
 an} (\nabla^F,g^{TM} )$ is a holomorphic section in the sense of
 \cite[Definition 3.4]{BK3}.

 By \cite[Theorem 11.3]{BK2} (cf. \cite[(5.13)]{BK3}),
 (\ref{3.13}), (\ref{3.161})  and (\ref{3.17}), one gets
 the following formula for the Ray-Singer norm of ${\cal T}_{\rm
 an} (\nabla^F,g^{TM} )$,
\begin{align}\label{3.18}
 \left\|{\cal T}_{\rm
 an}\left(\nabla^F,g^{TM}\right)\right\|^{\rm RS}=1.
  \end{align} In particular, when restricted to the space of acyclic representations,  ${\cal T}_{\rm
 an} (\nabla^F,g^{TM} )$ becomes a (univalent) holomorphic function such that
\begin{align}\label{3.19}
 \left|{\cal T}_{\rm
 an}\left(\nabla^F,g^{TM}\right)\right| =T^{RS}(\nabla^F),
  \end{align} the usual Ray-Singer analytic torsion. This provides an
  analytic resolution of a question of Burghelea (cf. \cite{BuH}).

  \begin{Rem} It is clear from the definition that ${\cal T}^2_{\rm
  an}$ does not depend on the choice of   $\widetilde{\nabla}^F$,
  and thus gives an intrinsic definition of a
 holomorphic section of the square of the determinant line bundle.
 \comment{The dependence of ${\cal T}_{\rm
 an}$ on $\alpha$ indicates in part the subtleness of the analytic meaning of the
 phase of the Turaev  torsion (cf. \cite{T, FT2}).}
 \end{Rem}

Next, we show how to   modify the Turaev torsion (cf. \cite{T,
FT2}) to get a holomorphic section with Ray-Singer norm equal to
one.

Let $\varepsilon$ be an Euler structure on $M$ and ${\bf o}$ a
cohomological orientation. We use the notation as in \cite{BK3} to
denote the associated Turaev torsion by $\rho_{\varepsilon,{\bf
o}}$.

Let $c(\varepsilon)\in H_1(M,{\bf Z})$ be the canonical class
associated to the Euler structure $\varepsilon$ (cf. \cite{T} or
\cite[Section 5.2]{FT2}). Then for any representation $\alpha_F$
corresponding to a flat vector bundle $(F,\nabla^F)$, by
\cite[Theorem 10.2]{FT2} one has
\begin{align}\label{3.20}
 \left\|\rho_{\varepsilon,{\bf
o}}\left(\alpha_F\right) \right\|^{\rm RS} =\left|{\rm
det}\,\alpha_F(c(\varepsilon))\right|^{1/2}.
  \end{align}

  Let ${\rm  L}_{\dim M-1}(TM)\in H^{\dim M-1}(M,{\bf Z})$ be the degree $\dim M-1$ component
  of the characteristic class ${\rm  L}(TM)$. Let $\widehat{\rm
  L}_1(TM)\in H_1(M,{\bf Z})$ denote its Poincar\'e dual. Then one
  verifies easily that
\begin{align}\label{3.21}
 \left|{\rm det}\, \alpha_F\left(\widehat{\rm
  L}_1(TM)\right)\right|=\exp\left(\int_M {\rm
L}\left(TM,\nabla^{TM}\right) c_1\left(F,\nabla^F\right)\right).
  \end{align}

  On the other hand, by
  \cite[Corollary 5.9]{BK3}, $\widehat{\rm
  L}_1(TM)+c(\varepsilon)\in H_1(M,{\bf Z})$ is divisible by two,
  and one can define a class $\beta_\varepsilon\in H_1(M,{\bf Z})$
  such that
\begin{align}\label{3.22}
  -2\beta_\varepsilon=\widehat{\rm
  L}_1(TM)+c(\varepsilon).
  \end{align}

  From Proposition \ref{t3.1}, (\ref{3.21}) and (\ref{3.22}), one finds
\begin{align}\label{3.23}
\left|{\rm det}\,\alpha_F(c(\varepsilon))\right|^{1/2}=\left|{\rm
det}\,\alpha_F\left(\beta_\varepsilon\right)\right|^{-1}\exp\left(-\pi\Phi\left(F,\nabla^F\right)
+\pi{\rm Im}\left(\overline{\eta}\left(D_{\rm
Sig}^F\right)\right)\right),
  \end{align}
  where $ \Phi(F,\nabla^F) $ is the locally constant function given
  by (\ref{3.12}).

  We now define a modified  Turaev torsion as follows:
\begin{align}\label{3.24}
 {\cal T}_{\varepsilon,{\bf o}}\left(F,\nabla^F\right)=\rho_{\varepsilon,{\bf
o}}\left(\alpha_F\right)e^{\pi\Phi\left(F,\nabla^F\right)+
\sqrt{-1}\pi\widetilde{\eta}\left(\nabla^F,\widetilde{\nabla}^F\right)}\left({\rm
det}\,\alpha_F\left(\beta_\varepsilon\right)\right) .
  \end{align}
  Clearly, ${\cal T}_{\varepsilon,{\bf o}} (F,\nabla^F )$ is a
  holomorphic section in the sense of \cite[Definition 3.4]{BK3}.
  Moreover, by (\ref{3.20}), (\ref{3.23}) and (\ref{3.24}), its
  Ray-Singer norm equals to one. Thus it provides  another
  resolution of Burghelea's problem mentioned above which should be closely related to what in \cite{BuH0}.

  Combining with (\ref{3.18}) we
  get
\begin{align}\label{3.25}
 \left| {{\cal T}_{\rm an}\left(\nabla^F, g^{TM}\right) \over
  {\cal T}_{\varepsilon,{\bf o}}\left(F,\nabla^F\right)}\right|=1,
    \end{align}
    which, in view of (\ref{3.11}), is equivalent to
    \cite[(5.10)]{BK3}.

    On the other hand, since now $ {{\cal T}_{\rm an} (\nabla^F, g^{TM} )
    /  {\cal T}_{\varepsilon,{\bf o}} (F,\nabla^F )}$ is a
    holomorphic function with absolute value identically equals to
    one, one sees that there is a real locally constant function
    $\theta_{\varepsilon,{\bf o}}(F,\nabla^F    )$ such that
\begin{align}\label{3.26}
  {{\cal T}_{\rm an}\left(\nabla^F, g^{TM}\right) \over
  {\cal T}_{\varepsilon,{\bf o}}\left(F,\nabla^F\right)}
   =e^{\sqrt{-1}\theta_{\varepsilon,{\bf o}}(F,\nabla^F    )},
    \end{align}
which is equivalent to \cite[(5.8)]{BK3}.

$\ $

\begin{Rem} While the univalent holomorphic sections ${\cal
T}_{\rm an}$ and ${\cal T}_{\varepsilon,{\bf o}}$ depend on the
choice of an ``initial'' flat connection  $\widetilde{\nabla}^F$,
the quotients in the left hand sides of (\ref{3.25}) and
(\ref{3.26}) do not involve it.
\end{Rem}

$\ $

\begin{Rem} One of the advantages of (\ref{3.25}) and (\ref{3.26})
is that they look in closer resemblance to the theorems of
Cheeger, M\"uller and Bismut-Zhang  (cf. \cite{BZ}) concerning the
Ray-Singer and Reidemeister torsions.
\end{Rem}

$\ $

Now let $\nabla^F_1$ and $\nabla^F_2$ be two acyclic unitary flat
connections on $F$. We do not assume that they can be connected by
a smooth path of flat connections.

By \cite[(14.11)]{BK1} (cf. \cite[(6.2)]{BK3}), (\ref{3.16}),
(\ref{3.17}) and the variation formula for $\eta$-invariants (cf.
\cite{APS1, APS3, BF}), one finds
\begin{align}\label{3.27}
  \frac{{\cal T}_{\rm an}\left(\nabla^F_1, g^{TM}\right)}
 {{\cal T}_{\rm an}\left(\nabla^F_2, g^{TM}\right)}
   = {T^{\rm RS}\left(\nabla^F_1\right)\over T^{\rm
   RS}\left(\nabla^F_2\right)}\cdot
   {\exp\left(-\sqrt{-1}\pi\overline{\eta}\left( D_{\rm
Sig,1}^F\right)+\sqrt{-1}\pi\widetilde{\eta}\left(\nabla^F_1,\widetilde{\nabla}^F\right)\right)\over
\exp\left(-\sqrt{-1}\pi \overline{\eta}\left(D_{\rm
Sig,2}^F\right)+\sqrt{-1}\pi\widetilde{\eta}\left(\nabla^F_2,\widetilde{\nabla}^F\right)\right)}
    \end{align}
    $$={T^{\rm RS}\left(\nabla^F_1\right)\over T^{\rm
   RS}\left(\nabla^F_2\right)}\cdot {\exp\left(-\sqrt{-1}\pi\overline{\eta}\left( D_{\rm
Sig,1}^F\right)+\sqrt{-1}\pi \overline{\eta}\left(D_{\rm
Sig,2}^F\right)\right)\over \exp\left(-\sqrt{-1}\pi\int_M{\rm
L}\left(TM,\nabla^{TM}\right)CS\left( {\nabla}^{F
}_2,\nabla^F_1\right)\right)}$$
$$={T^{\rm RS}\left(\nabla^F_1\right)\over T^{\rm
   RS}\left(\nabla^F_2\right)}\cdot\exp\left(\sqrt{-1}\pi\cdot {\rm sf}\left( D_{\rm
Sig,1}^F, D_{\rm Sig,2}^F\right)\right),$$ where $D_{\rm Sig,1}^F$
and $ D_{\rm Sig,2}^F$ are the signature operators associated to
$\nabla^F_1$ and $\nabla^F_2$ respectively, while ${\rm sf} (
D_{\rm Sig,1}^F, D_{\rm Sig,2}^F )$ is the spectral flow of the
linear path connecting $D_{\rm Sig,1}^F$ and  $D_{\rm Sig,2}^F$,
in the sense of Atiyah-Patodi-Singer \cite{APS3}.

\begin{Rem}
Since we do not assume that $\nabla^F_1$ and $\nabla^F_2$ can be
connected by a path of flat connections, our formula extends the
corresponding formula in \cite[Proposition 6.2]{BK3}.
\end{Rem}

\begin{cor} The ratio ${{\cal T}_{\rm an} (\nabla^F,g^{TM} )/ T^{\rm
   RS} (\nabla^F )}$ is a locally constant function on the set of
   acyclic unitary flat connections on $F$.
\end{cor}

\begin{ex} Let $\nabla^F$ be an acyclic unitary flat connection on $F$.
Let $g\in \Gamma( U(F))$ be a smooth section of unitary
automorphisms of $F$. Then $g^{-1} \nabla^Fg$ is another acyclic
unitary flat  connection on $F$. A  standard calculation shows
that
\begin{align}\label{3.28}
 {\rm sf}\left(D_{\rm
Sig}^{F,\nabla^F}, D_{\rm Sig}^{F,g^{-1}\nabla^Fg}\right)=\int_M
{\rm L}(TM){\rm ch}(g),
    \end{align} where ${\rm ch}(g)\in H^{\rm odd}(M,{\bf R})$ is the odd Chern character
    associated to $g$ (cf. \cite{Z1}).  From
    (\ref{3.28}), one sees that if
    $\int_M
{\rm L}(TM){\rm ch}(g)$ is nonzero, then $\nabla^F$ and
$g^{-1}\nabla^Fg$ do not lie in the same connected  component in
the set of
   acyclic unitary flat connections on $F$.
\end{ex}

\subsection{More on $\eta$-invariants, spectral flow and the phase of the refined analytic torsion}

We would like to point out that the (reduced) $\eta$-invariant for
non-self-adjoint operators  we used above, when considered as a
${\bf C} $-valued function, is the original $\eta$ invariant
appeared in \cite{APS3} (see also \cite{G}). In this section, we
show that the ${\bf R}$-valued variation formula for
$\eta$-invariants (which has been used in (\ref{3.27})) admits an
extension to a ${\bf C}$-valued variation formula valid also for
the non-self-adjoint operators discussed in the present paper.

First, the concept of spectral flow can be extended to
non-self-adjoint operators, and this has been done in \cite{ZL} in
a general context.

For our specific situation, if $D_{{\rm Sig}, t}^{F}$, $0\leq
t\leq 1$, is a smooth path of (possibly) non-self-adjoint
signature operators, following \cite{ZL}, we define the spectral
flow of this path to be, tautologically,
\begin{align}\label{3.29}
 {\rm sf}\left(D_{\rm Sig, 0}^{F},D_{\rm Sig,
 1}^{F}\right)=\#\left\{\mbox{${\rm spec}\left(D_{\rm Sig,
 0}^{F }\right)\cap\left\{{\rm Re}(\lambda)\geq 0\right\} $ $\rightarrow$ ${\rm spec}\left(D_{\rm Sig,
 1}^{F }\right)\cap\left\{{\rm Re}(\mu)<0\right\}$}\right\}
    \end{align}
    $$-\#\left\{\mbox{${\rm spec}\left(D_{\rm Sig,
 0}^{F }\right)\cap\left\{{\rm Re}(\lambda)< 0\right\} $ $\rightarrow$ ${\rm spec}\left(D_{\rm Sig,
 1}^{F }\right)\cap\left\{{\rm Re}(\mu)\geq 0\right\}$}\right\},$$
 which simply replaces the number zero in 
the original definition for self-adjoint
 operators (\cite{APS3}) by the  axis of purely imaginary numbers.

 Now let $\nabla^F_t$, $0\leq t\leq 1$, be a smooth path of (not
 necessary unitary and/or flat)
 connections on $F$. Let $D_{{\rm Sig}, t}^{F}$, $0\leq
t\leq 1$, be the corresponding path of signature operators. With
the definition of spectral flow, one then sees easily that the
following
 variation formula holds in ${\bf C}$,
\begin{align}\label{3.30}
\overline{\eta}\left( D_{\rm Sig,
 1}^{F }\right)-\overline{\eta}\left( D_{\rm Sig,
 0}^{F }\right)={\rm sf}\left(D_{\rm Sig, 0}^{F},D_{\rm Sig,
 1}^{F}\right)+\int_M{\rm
 L}\left(TM,\nabla^{TM}\right)CS\left(\nabla^F_0,\nabla^F_1\right).
    \end{align}

    Now we observe that in \cite{BK1, BK2, BK3}, Braverman and
    Kappeler propose an alternate definition of (reduced) $\eta$
    invariant, which if we denote by $\eta_{BK}$, then (cf.
    \cite[Definition 4.3]{BK1} and \cite[Definition 5.2]{BK3})
\begin{align}\label{3.31}
\eta_{BK}\left( D_{\rm Sig}^{F }\right)=\overline{\eta}\left(
D_{\rm Sig}^{F }\right)-m_{-}\left(D_{\rm Sig}^{F }\right),
    \end{align}
    where $m_{-} (D_{\rm Sig}^{F } )$ is the number of purely imaginary eigenvalues of
    $D_{\rm Sig}^{F }$ of form $\lambda\sqrt{-1}$ with $\lambda<0$.

Formulas (\ref{3.30}) and (\ref{3.31}) together give a variation
formula for $\eta_{BK}$, which can be used to extend (\ref{3.27})
to non-unitary acyclic representations.

\end{document}